\documentclass[12pt]{article}

\usepackage[english]{babel}

\usepackage{amsfonts, amssymb, amsmath, amsthm, epsf, cite}

\usepackage{graphicx}

\usepackage{fullpage}

\newtheorem{theorem}{Theorem}[section]
\newtheorem{lemma}{Lemma}[section]

\theoremstyle{definition}
\newtheorem{definition}{Definition}[section]

\newtheorem{remark}{Remark}[section]

\numberwithin{equation}{section} \numberwithin{figure}{section}

\newcommand{\ph}{{\varphi}}

\newcommand{\cN}{{\mathcal N}}

\newcommand{\cH}{{H}}
\newcommand{\cF}{{F}}

\newcommand{\bdelta}{{\delta}}
\newcommand{\bDelta}{{\Delta}}

\newcommand{\bu}{u}
\newcommand{\br}{{\mathbf r}}
\newcommand{\bv}{v}

\newcommand{\bbR}{{\mathbb R}}
\newcommand{\bbN}{{\mathbb N}}

\newcommand{\bbZ}{{\mathbb Z}}

\let\phi=\varphi

\newcommand{\ep}{\varepsilon}

\title{Asymptotics of the heat kernels on 2D lattices}

\author{Pavel  Gurevich\footnote{Free University of Berlin, Germany; Peoples' Friendship
University of Russia, Russia; email: gurevich@math.fu-berlin.de}}

\begin{document}

\maketitle

\begin{abstract}
We obtain asymptotic expansions of the spatially discrete 2D heat kernels, or Green's functions on lattices, with respect to powers of time variable up to an arbitrary order and estimate the remainders uniformly
on the whole lattice. Unlike in the 1D case, the asymptotics contains a time independent term. The derivation of its spatial asymptotics is the technical core of the paper. Besides numerical applications, the obtained results play a crucial role in the analysis of spatio-temporal patterns for reaction-diffusion equations on lattices, in particular  rattling patterns for hysteretic diffusion systems.
\end{abstract}

{\bf Keywords:} discrete heat kernel, Green's function, lattice dynamics, asymptotics

{\bf MSC:} 39A70, 35K08, 35B40, 34A33 

\tableofcontents

\section{Introduction}

The paper deals with the spatially discrete heat kernels, or Green's functions of the heat equation, on 2D lattices. Although one can easily obtain integral representations of the discrete Green's function, using the discrete Fourier transform, it is not
possible to express them via elementary functions. Therefore, their
asymptotic expansions play an important role in applications. Asymptotics of 
 lattice
 Green's functions in the stationary (elliptic) case were studied beginning from 1950s, see~\cite{Duffin}, the subsequent papers~\cite{Bramble, Katsura,
Mangad, Martinsson, Guttmann, Otto}, and the monograph~\cite[Chapter~8]{Lifanov}.
For parabolic operators, there is vast literature in the spatially continuous case. For example, large-time behavior of Green's functions was treated in~\cite{Murata85} (for small perturbations of the heat operator) and in \cite{Tsuchida08,Norris97} (for spatially periodic coefficients). A survey on the large time behavior of heat kernels for second-order parabolic operators on Riemannian manifolds can be found in~\cite{Pinchover13}. In the spatially discrete case, the research directions include
continuous-time random walks on general graphs (see, e.g.,~\cite{Lawler,
Pang93} and references therein) and on lattices in a random
environment (see, e.g.,~\cite{Delmotte}). In both cases,
Gaussian bounds for the heat kernel are extensively studied. However, higher-order asymptotics of  Green's functions are not available in general. We mention~\cite{Iliev02, Iliev07}, where an asymptotic expansion of  Green's function for particular parabolic equations on 1D lattices was obtained in terms of the Bessel functions. In the case of constant coefficients, we obtained asymptotic formulas as $t\to\infty$ for Green's functions of general higher-order parabolic operators, with uniform estimates of the remainders on the whole lattices~\cite{GurAsympGeneral}. These formulas played a crucial role in explaining spatio-temporal patterns (rattling) for hysteretic reaction-diffusion equations on 1D lattices~\cite{GurTikhRattlingMMS,GurTikhRattlingPoincare}. To generalize those results to 2D lattices, one needs explicit asymptotic expansions of 2D heat kernels, which is one of the main motivations for the present paper. Unlike in the 1D case, the general asymptotics in~\cite{GurAsympGeneral} contains a time independent term $\Omega(x/\varepsilon)$ given in the integral form, where $\ep$ is the grid step. Derivation of its spatial asymptotics is a technical core of our paper. It turns out that, unlike the leading time-dependent term, $\Omega(x/\varepsilon)$ is not rotationally symmetric but rather depends on the polar angle of $x$. Furthermore, it vanishes  for $x\ne0$ as $\varepsilon\to 0$ and, as such, is an artifact of spatial discretization (generally depending on the lattice structure).

The paper is organized as follows. In Sec.~\ref{secNotation}, we introduce general notation and define first and second Green's functions on 2D lattices. In Sec.~\ref{secAsympFirstGreenFunction}, we obtain a theorem on asymptotics of the first Green's function, which directly follows from a general result in~\cite{GurAsympGeneral}. In Sec.~\ref{secAsympSecondGreenFunction}, we formulate a theorem on asymptotics of the second Green's function, which contains the above-mentioned term $\Omega(x/\varepsilon)$. Derivation of asymptotics of the latter is contained in Sec.~\ref{secProof}. Appendix~\ref{appendixAuxiliary} contains auxiliary results on the Bessel functions of the first kind, which are used in the main part of the paper.

\section{Notation}\label{secNotation}

Green's functions are special solutions of the heat equations  on the 2D {\it grid space}, or {\it
lattice},
$$
\bbR_\ep^2:=\{x\in\bbR^2: x=\ep s,\ s\in\bbZ^2\},\ \ep>0.
$$
The heat equations are defined via the Laplace operator on $\bbR_\ep^2$ given by
$$
\bDelta_\ep w(x):=\ep^{-2}\big(w(x_1+\ep,x_2)+w(x_1-\ep,x_2) + w(x_1,x_2+\ep)+ w(x_1,x_2-\ep)-4w(x)\big),\quad x\in\bbR_\ep^2.
$$
Let $\bdelta^\ep(x)$ denote the {\it
grid delta-function} given by
\begin{equation}\label{eqdelta}
\bdelta^\ep(0)=\ep^{-2},\qquad \bdelta^\ep(x)=0\quad \forall
x\in\bbR_\ep^2\setminus\{0\}.
\end{equation}

\begin{definition}\label{defGreenFirst}
  We call the function $\bu_\ep(x,t)$, $x\in\bbR_\ep^2,\ t\ge 0$,
the {\em first discrete Green function} if $\bu_\ep(\cdot,t)$ is a
rapidly decreasing grid function for all $t\ge0$, $\bu_\ep(x,\cdot)\in
C^1[0,\infty)$ for all $x\in\bbR_\ep^2$, and
$$
\left\{
\begin{aligned}
& \dot\bu_\ep(x,t)-\bDelta_\ep \bu_\ep(x,t)=0, & & x\in\bbR_\ep^2,\ t>0,\\
& \bu_\ep(x,0)=\bdelta_\ep(x),   & & x\in\bbR_\ep^2.
\end{aligned}
\right.
$$
\end{definition}

\begin{definition}\label{defGreenSecond}
We call the function $\bv_\ep(x,t)$, $x\in\bbR_\ep^2,\ t\ge 0$, the {\em
second discrete Green function} if $\bv_\ep(\cdot,t)$ is a rapidly
decreasing grid function for all $t\ge0$, $\bv_\ep(x,\cdot)\in
C^1[0,\infty)$ for all $x\in\bbR_\ep^2$, and
$$
\left\{
\begin{aligned}
& \dot\bv(x,t)-\bDelta_\ep \bv(x,t)=\bdelta_\ep(x), & & x\in\bbR_\ep^2,\ t>0,\\
& \bv_\ep(x,0)=0,   & & x\in\bbR_\ep^2.
\end{aligned}
\right.
$$
\end{definition}

Using the discrete Fourier transform (see~\cite{GurAsympGeneral} for details), we obtain the explicit
representations
\begin{align}
\bu_\ep(x,t)&=\dot \bv_\ep(x,t)
=\dfrac{1}{(2\pi)^2}\int_{R_{\pi\ep^{-1}}}
e^{-t\ep^{-2} A(\eta\ep)}e^{i
x\eta}\,d\eta,\label{eqz_nFourierMultiDim}\\
\bv_\ep(x,t)&=\dfrac{1}{(2\pi)^2}\int_{R_{\pi\ep^{-1}}}
\dfrac{1-e^{-t\ep^{-2} A(\eta\ep)}}{\ep^{-2}A(\eta\ep)}\,e^{i
x\eta}\,d\eta,\label{eqy_nFourierMultiDim}
\end{align}
where
\begin{equation}\label{eqSymbolA}
A(\xi):=2(2-\cos\xi_1-\cos\xi_2)
\end{equation}
is the {\em symbol} of the operator $-\bDelta_1$ and
$$
R_{L}:=\left\{\eta\in\bbR^2: |\eta_k|\le L,\
k=1,2\right\}.
$$

\begin{remark}\label{remEp1}
  Changing the variables in the integrals
in~\eqref{eqz_nFourierMultiDim} and~\eqref{eqy_nFourierMultiDim},
we obtain for $J=0,1,2,\dots$
\begin{equation}\label{eqRescalingGreen}
\dfrac{\partial ^J \bu^\ep(x,t)}{\partial t^J}\equiv
\ep^{-2(J+1)} \dfrac{\partial^J
\bu^1\left(x',\tau\right)}{\partial
\tau^J}\Bigg|_{x'=x/\ep,\,\tau=t/\ep^2},\quad
\bv^\ep(x,t)\equiv
\bv^1\left(\dfrac{x}{\ep},\dfrac{t}{\ep^2}\right).
\end{equation}
\end{remark}

\section{Asymptotics of first Green's function}\label{secAsympFirstGreenFunction}

First, we formulate a theorem on asymptotics of the first Green's function $\bu_\ep(x)$.

\begin{theorem}\label{thAsymphMultiDim}
For any $\ep>0$, $t_0>0$, integer $N\ge 1$, integer $J\ge 0$, and
all $x\in\bbR_\ep^2$ and $t\ge t_0\ep^2$, we have
\begin{equation}\label{eqAsymphGeneralMultiDim}
\dfrac{\partial^J\bu_\ep(x,t)}{\partial t^J}=\sum\limits_{n=0}^{N-1}
\dfrac{\ep^{2n}}{t^{n+1+J}}\,
\cH_{Jn}\left(\dfrac{x}{t^{1/2}}\right) + \br_\bu^\ep(J,N,t_0;
x,t),
\end{equation}
where
\begin{equation}\label{eqHk}
\cH_{Jn}(y):=\left\{
\begin{aligned}
&\Delta^J\cH(y) & &\text{if } n=0,\\
&\text{finite linear combinations}\\
&\text{of derivatives of $\cH(y)$} & &\text{if } n=1,\dots,N-1,
\end{aligned}
\right.
\end{equation}
$$\Delta=\frac{\partial^2}{\partial y_1^2}+\frac{\partial^2}{\partial y_2^2},$$
\begin{equation}\label{eqhGeneral}
\cH(y)=\dfrac{1}{4\pi}
e^{-|y|^2/4},
\end{equation}
\begin{equation}\label{eqHkRemainder}
|\br_\bu^\ep(J,N,t_0; x,t)|\le
\dfrac{\ep^{2N}R_\bu(J,N,t_0)}{t^{N+1+J}},
\end{equation}
and $R_\bu(J,N,t_0)\ge 0$ does not depend on
$x\in\bbR_\ep^2$, $t\ge t_0\ep^2$, and $\ep>0$.
\end{theorem}
\proof The result follows from~\cite[Theorem~3.1]{GurAsympGeneral}. Note that the corresponding asymptotic expansion in~\cite[Theorem~3.1]{GurAsympGeneral} is obtained for parabolic problems with general higher-order elliptic parts and additionally contains fractional powers $t^{k/2+1+J}$ with odd $k$. The corresponding factors $H_{J,k/2}(\cdot)$ correspond to homogeneous polynomials of odd degrees in the Taylor expansion of the symbol $A(\xi)$ of the elliptic part. In our case, these odd degree polynomials are absent, see~\eqref{eqSymbolA}. Hence, the corresponding factors $H_{J,k/2}(\cdot)$ in the asymptotic expansion vanish too.
\endproof

\begin{remark}
\begin{enumerate}
\item The first term in the asymptotics of $\bu_\ep(x,t)$ is $\dfrac{1}{4\pi t}
e^{-|x|^2/(4t)}$, which coincides with Green's function of the continuous heat operator.

\item All the functions $H_{Jn}(y)$ entering the higher order terms can be found explicitly, see~\cite[Remark~3.1]{GurAsympGeneral} for details. For example,
$$
\begin{aligned}
H_{01}(y)&=\frac{2}{4!}\left(\frac{\partial^4}{\partial y_1^4}+\frac{\partial^4}{\partial y_2^4}\right)H(y),\\ H_{02}(y)&=\frac{2}{6!}\left(\frac{\partial^6}{\partial y_1^6}+\frac{\partial^6}{\partial y_2^6}\right)H(y) + \frac{4}{2!(4!)^2}\left(\frac{\partial^4}{\partial y_1^4}+\frac{\partial^4}{\partial y_2^4}\right)^2H(y).
\end{aligned}
$$
\end{enumerate}
\end{remark} 

\section{Asymptotics of second Green's function}\label{secAsympSecondGreenFunction}

In this section, we formulate a theorem on asymptotics of the second Green's function $\bv_\ep(x,t)$. The proof of this theorem is given in Sec.~\ref{secProof}.

For $y\in\bbR^2\setminus\{0\}$, we introduce the functions
\begin{equation}\label{eqcFLog-Corrected2}
\begin{aligned}
\cF_0(y) & :=
 2   \int_r^\infty
\dfrac{\cH_{00}(\rho,\ph)}{\rho}\,d\rho=\frac{1}{2\pi}\int_{r}^{\infty}\frac{e^{-\rho^2/4}}{\rho} d\rho,\\
\cF_n(y) & := -\dfrac{2}{r^{2n}}  \int_0^r
  \rho^{2n-1}\cH_{0n}(\rho,\ph)\,d\rho,\quad
n\ge 1.
\end{aligned}
\end{equation}
Here and below we denote functions written in Cartesian and spherical
coordinates by the same letter, e.g., $\cH_{0n}(r,\ph)$ stands
for $\cH_{0n}(y)$.

In what follows, we set $B_1:=\{\theta\in\bbR^2: |\theta|<1\}$. We also introduce the function
\begin{equation}\label{eqr_pi}
r_\pi(\ph):=\text{distance from the origin to the boundary of
 the square $R_\pi$ in the direction $\ph$}.
\end{equation}

\begin{theorem}\label{thAsympfMultiDim}
For any $\ep>0$, $t_0>0$, integer $N\ge 1$, and
 $t\ge t_0\ep^2$, the following holds.
\begin{enumerate}

\item\label{thAsympfMultiDim1} If $x\in\bbR_\ep^2\setminus\{0\}$, then
\begin{equation}\label{eqAsympfGeneralMultiDim}
\bv_\ep(x,t)=\cF_{0}\left(\dfrac{x}{t^{1/2}}\right) + \Omega\left(\frac{x}{\ep}\right)+\sum\limits_{n=1}^{N-1}
\dfrac{\ep^{2n}}{t^{n}}\,
\cF_{n}\left(\dfrac{x}{t^{1/2}}\right) + \br_\bv^\ep(N,t_0;
x,t),
\end{equation}
where $\cF_n(y)$ are given by~\eqref{eqcFLog-Corrected2},
\begin{equation}\label{eqOmega1Asymp}
\Omega\left(\frac{x}{\ep}\right) =\dfrac{\ep^2}{24\,\pi}\dfrac{\cos(4\psi)}{r^2} + \br_\Omega^\ep(x),
\end{equation}
$(r,\psi)$ are the polar coordinates of $x$,
\begin{equation}\label{eqOmegaRemainder}
|\br_\Omega^\ep(x)|\le
\dfrac{\ep^{5/2} R_\Omega}{r^{5/2}},
\end{equation}
and $R_\Omega\ge0$ does not depend on $x\in\bbR_\ep^2\setminus\{0\}$ and $\ep>0$.

\item\label{thAsympfMultiDim2} If $x=0$, then
$$
\bv(0,t)=\cH(0)\ln\frac{t}{\ep^2} + S_0 - \sum\limits_{n=1}^{N-1}
\dfrac{\ep^{2n}}{t^{n}}\,
\dfrac{H_{0n}(0)}{n} + \br_\bv^\ep(N,t_0;0,t),
$$
where $H(\cdot)=H_{00}(\cdot)$ and $H_{0n}(\cdot)$ are the functions in~\eqref{eqhGeneral} and~\eqref{eqHk},
\begin{equation}\label{eqS0}
S_0=\dfrac{1}{(2\pi)^2}  \left(\pi\gamma+ \int\limits_{B_1}\left(\dfrac{1}{A(\theta)}-\dfrac{1}{|\theta|^2}\right) d\theta + \int_0^{2\pi} \ln r_\pi(\ph)\,d\ph\right),
\end{equation}
$r_\pi(\ph)$ is given by~\eqref{eqr_pi}, and $\gamma$ is the Euler--Mascheroni constant.
\end{enumerate}
In both cases,
\begin{equation}\label{eqFkRemainder}
|\br_\bv^\ep(N,t_0; x,t)|\le
\dfrac{\ep^{2N}R_\bv(N,t_0)}{t^{N}}
\end{equation}
and  $R_\bv(N,t_0)\ge 0$ does not depend on
$x\in\bbR_\ep^2$, $t\ge t_0$, and $\ep>0$.
\end{theorem}

Item~\ref{thAsympfMultiDim2} in Theorem~\ref{thAsympfMultiDim1} follows from~\cite[Theorem~5.2, part 2]{GurAsympGeneral}, in which the constant $S_0$ is given by
\begin{equation}\label{eqS01}
\begin{aligned}
S_0  :=\dfrac{1}{(2\pi)^2}  \Bigg( &
 \int\limits_{B_1}\dfrac{1-e^{-|\xi|^2}}{|\xi|^2}\,d\xi -
\int\limits_{\bbR^2\setminus
 B_1}\dfrac{e^{-|\xi|^2}}{|\xi|^2}\,d\xi
\\
& + \int\limits_{R_\pi}\left(\dfrac{1}{A(\theta)}-\dfrac{1}{|\theta|^2}\right) d\theta +  \int\limits_0^{2\pi} \ln r_\pi(\ph)\,d\ph   \Bigg).
\end{aligned}
\end{equation}
Passing to the polar coordinates and using Lemma~\ref{lExpGamma}, we see that
\begin{equation}\label{eqS02}
\begin{aligned}
  \int\limits_{B_1}\dfrac{1-e^{-|\xi|^2}}{|\xi|^2}\,d\xi -
 \int\limits_{\bbR^2\setminus
 B_1}\dfrac{e^{-|\xi|^2}}{|\xi|^2}\,d\xi & =2\pi\left( \int_0^1\dfrac{1-e^{-\rho^2}}{\rho}\,d\rho -
 \int_1^\infty\dfrac{e^{-\rho^2}}{\rho}\,d\rho\right)\\
 & =\pi\left( \int_0^1\dfrac{1-e^{-z}}{z}\,dz -
 \int_1^\infty\dfrac{e^{-z}}{z}\,dz\right)=\pi \gamma,
 \end{aligned}
\end{equation}
where $\gamma$ is the Euler--Mascheroni constant. Combining~\eqref{eqS01} and~\eqref{eqS02} yields~\eqref{eqS0}.

The proof of item~\ref{thAsympfMultiDim1} is given in Sec.~\ref{secProof}.

\section{Proof of Theorem~$\ref{thAsympfMultiDim}$}\label{secProof}

\subsection{Integral representation of $\Omega(x)$}

Taking into account Remark~\ref{remEp1}, it suffices to prove item~\ref{thAsympfMultiDim1} in Theorem~\ref{thAsympfMultiDim1} for $\ep=1$. Therefore, from now on, we consider $x\in\bbZ^2$.

Due to~\cite[Theorem~5.2, part 1]{GurAsympGeneral}, we have
\begin{equation}\label{eqOmega1}
\begin{aligned}
\Omega(x) =\dfrac{1}{(2\pi)^2}  \Bigg(& 2 \int\limits_{B_1}\dfrac{1-e^{-|\xi|^2}}{|\xi|^2}\,d\xi -
2 \int\limits_{\bbR^2\setminus
 B_1}\dfrac{e^{-|\xi|^2}}{|\xi|^2}\,d\xi + \int\limits_{R_\pi}\left(\dfrac{e^{i x\theta}}{A(\theta)}-\dfrac{1}{|\theta|^2}\right) d\theta \\
 & +  \int\limits_0^{2\pi} \ln(r_\pi(\ph))\,d\ph   \Bigg) +\dfrac{1}{2\pi}(\ln r - \ln 2).
\end{aligned}
\end{equation}
Hence, our main goal is to prove~\eqref{eqOmega1Asymp}, \eqref{eqOmegaRemainder} for $\Omega(x)$ given by~\eqref{eqOmega1}.

We begin with simplifying~\eqref{eqOmega1}. Writing $x\in\bbR^2\setminus\{0\}$ in the polar coordinates $(r,\psi)$ and $\theta\in\bbR^2$ in the polar coordinates $(\rho,\phi)$, we have
\begin{equation}\label{eqOmega3}
\int\limits_{R_\pi}\left(\dfrac{e^{i x\theta}}{A(\theta)}-\dfrac{1}{|\theta|^2}\right) d\theta=I_1(r,\psi)+I_2(r,\psi),
\end{equation}
where
\begin{equation}\label{eqI1}
I_1(r,\psi):=
\int_{-\pi}^{\pi}d\phi \int_0^{r_\pi(\phi)}
\cos(r\rho\cos(\phi-\psi))\left(\dfrac{1}{A(\rho\cos\phi,\rho\sin\phi)}-\dfrac{1}{\rho^2}\right)\rho\,d\rho,
\end{equation}
$$
I_2(r,\psi):=
\int_{-\pi}^{\pi}d\phi \int_0^{r_\pi(\phi)}
\dfrac{\cos(r\rho\cos(\phi-\psi))-1}{\rho}\,d\rho.
$$
In the integral defining $I_2$, we change the variable $z=r\rho$ and   obtain
$$
\begin{aligned}
I_2(r,\psi)&=\int_{-\pi}^{\pi}d\phi \int_0^{r_\pi(\phi)r}
\dfrac{\cos(z\cos(\phi-\psi))-1}{z}\,dz \\
&= \int_{-\pi}^{\pi}d\phi \int_0^1
\dfrac{\cos(z\cos(\phi-\psi))-1}{z}\,dz + \int_{-\pi}^{\pi}d\phi \int_1^{\pi r}
\dfrac{\cos(z\cos(\phi-\psi))}{z}\,dz \\
& + \int_{-\pi}^{\pi}d\phi \int_{\pi r}^{r_\pi(\phi)r}
\dfrac{\cos(z\cos(\phi-\psi))}{z}\,dz
- \int_{-\pi}^{\pi}d\phi \int_1^{r_\pi(\phi)r}
\dfrac{1}{z}\,dz.
\end{aligned}
$$
Using the integral representation of the Bessel function of the first kind $J_0(z)$ (see~\eqref{eqBesselIntegralRepresentation1}) and the identity from Lemma~\ref{lExpGamma}, we have
\begin{equation}\label{eqI2}
\begin{aligned}
& I_2(r,\psi)\\
&\quad = 2\pi  \int_0^1
\dfrac{J_0(z)-1}{z}\,dz + 2\pi \int_1^{\pi r}
\dfrac{J_0(z)}{z}\,dz + I_3(r,\psi) - \int_{-\pi}^{\pi} \ln(r_\pi(\phi))\,d\phi - 2\pi \ln r\\
& \quad =- 2\pi\gamma + 2\pi\ln 2 +I_4(r) + I_3(r,\psi) - \int_{-\pi}^{\pi} \ln(r_\pi(\phi))\,d\phi - 2\pi \ln r,
\end{aligned}
\end{equation}
where
\begin{equation}\label{eqI3}
I_3(r,\psi):=\int_{-\pi}^{\pi}d\phi \int_{\pi r}^{r_\pi(\phi)r}
\dfrac{\cos(z\cos(\phi-\psi))}{z}\,dz,\quad I_4(r):=-2\pi \int_{\pi r}^{\infty}
\dfrac{J_0(z)}{z}\,dz
\end{equation}

Combining~\eqref{eqOmega1}, \eqref{eqS02}, \eqref{eqOmega3}, \eqref{eqI2}, we obtain
\begin{equation}\label{eqOmegaI1I3I4}
\Omega(x)=\dfrac{1}{(2\pi)^2}  \big(I_1(r,\psi)+I_3(r,\psi)+I_4(r) \big).
\end{equation}

\subsection{Estimates of $I_1$, $I_3$, and $I_4$}

We recall that $(\rho,\phi)$ and $(r,\psi)$ stand for the polar coordinates of $\theta\in\bbR^2$ and $x\in\bbZ^2$, respectively. Set
\begin{equation}\label{eqI3:0}
\Psi_1:=[-\pi/4,\pi/4]\cup[3\pi/4,\pi]\cup[-\pi,-3\pi/4],\quad \Psi_2:=[0,2\pi)\setminus\Psi_1.
\end{equation}
Thus, the set $\{x\in\bbZ^2:\psi\in\Psi_j\}$ is the intersection of $\bbZ^2$ with the bisector of opening $\pi/2$ symmetric with respect to the axis $x_j$.

\begin{lemma}\label{lI3} If $\psi\in\Psi_1$, then
 \begin{equation}\label{eqI3'}
 I_3(r,\psi)= -\dfrac{2\sqrt{2}}{\pi} \sin\left(\pi r-\dfrac{\pi}{4}\right)  r^{-3/2} + \hat I_3(r,\psi) + O(r^{-5/2})\quad\text{as } r\to\infty,
 \end{equation}
   where $O(\cdot)$ is uniform with respect to $\psi\in[-\pi,\pi)$ and
 \begin{equation}\label{eqI3''}
 \hat I_3(r,\psi):= - \dfrac{1}{r \cos\psi}\int\limits_{R_\pi\setminus B_\pi} \sin(x\theta) \left(\dfrac{1}{|\theta|^2}\right)_{\theta_1}\,d\theta.
 \end{equation}
 If $\psi\in\Psi_2$, then  $I_3(r,\psi)= I_3(r,\psi+\pi/2)$.
\end{lemma}
\proof
 Assume that $\psi\in\Psi_1$   and hence $|\cos\psi|\ge 1/\sqrt{2}$. We rewrite $I_3(r,\psi)$ as follows:
\begin{equation}\label{eqI3:1'}
I_3(r,\psi)=\int\limits_{R_\pi\setminus B_\pi} \dfrac{\cos(x\theta)}{|\theta|^2}\,d\theta=\dfrac{1}{r\cos\psi} \int\limits_{R_\pi\setminus B_\pi} \dfrac{(\sin(x\theta))_{\theta_1}}{|\theta|^2}\,d\theta.
\end{equation}
Integrating by parts yields
\begin{equation}\label{eqI3:1}
I_3(r,\psi)=\tilde I_3(r,\psi)+\hat I_3(r,\psi),
\end{equation}
where
$$
\tilde I_3(r,\psi):=-\dfrac{1}{\pi r \cos\psi} \int_{-\pi}^\pi  \sin(\pi r\cos(\phi-\psi))\cos\phi\,d\phi, \\
$$
and $\hat I_3(r,\psi)$ is defined in~\eqref{eqI3''}.

Using the integral representation~\eqref{eqBesselIntegralRepresentation2} and Lemma~\ref{lAsympBesseln} with $n=1$, we obtain
\begin{equation}\label{eqI3:2}
\begin{aligned}
\tilde I_3(r,\psi)& =-\dfrac{1}{\pi r\cos\psi} \int_{-\pi}^\pi  \sin(\pi r\cos\phi)\cos\phi\cos\psi\,d\phi=-\dfrac{2}{r} J_1(\pi r)\\
&=-\dfrac{2\sqrt{2}}{\pi} \sin\left(\pi r-\dfrac{\pi}{4}\right)  r^{-3/2} + O(r^{-5/2}).
\end{aligned}
\end{equation}

To estimate $\hat I_3(r,\psi)$, we integrate by parts again and obtain $\hat I_3(r,\psi)=O(r^{-2})$. Together with~\eqref{eqI3:1} and~\eqref{eqI3:2}, this proves the lemma whenever $\psi\in\Psi_1$. If $\psi\in\Psi_2$, it suffices to note that $I_3(r,\psi)=I_3(r,\psi+\pi/2)$.
\endproof

\begin{lemma}\label{lI4}
  $I_4(r)= \dfrac{2\sqrt{2}}{\pi} \sin\left(\pi r-\dfrac{\pi}{4}\right)  r^{-3/2} + O(r^{-5/2})$ as $r\to\infty$.
\end{lemma}
\proof
The proof follows from the asymptotic expansion in Lemma~\ref{lAsympBessel0}.
\endproof

To write down an asymptotic formula for $I_1(r,\psi)$ (defined in~\eqref{eqI1}), we need the following notation for $\theta\ne 0$:
\begin{equation}\label{eqD0}
D_0(\theta):=\dfrac{1}{A(\theta)}-\dfrac{1}{|\theta|^2},
\end{equation}
\begin{equation}\label{eqD}
D(\theta):= \frac{\partial D_0(\theta)}{\partial\theta_1}=-\dfrac{2\sin\theta_1}{(4-2\cos\theta_1-2\cos\theta_2)^2}+\dfrac{2 \theta_1}{(\theta_1^2+\theta_2^2)^2},
\end{equation}
and
\begin{equation}\label{eqE}
E(\rho,\phi):=\rho D(\rho\cos\phi,\rho\sin\phi).
\end{equation}

Obviously, the functions $D_0(\theta)$ and $D(\theta)$ are infinitely differentiable for $\theta\in R_\pi\setminus\{0\}$, while $E(\rho,\phi)$ is infinitely differentiable for $\rho>0$ and $\phi\in\bbR$ and is $2\pi$ periodic with respect to $\phi$. It is also easy to see that $D_0(\theta)$ is bounded near the origin. The behavior of $E(\rho,\phi)$ near $\rho=0$ is described in the following lemma.

\begin{lemma}\label{lESmoothness}
 For any integer $k_1,k_2\ge 0$, the derivative $\dfrac{\partial^{k_1+k_2}E(\rho,\phi)}{\partial\rho^{k_1}\partial\phi^{k_2}}$ extends as a continuous function from $\{\rho>0,\phi\in\bbR\}$ to $\{\rho\ge 0,\phi\in\bbR\}$. Moreover,
 \begin{equation}\label{eqESmoothness}
  \lim\limits_{\rho\to 0}E(\rho,\phi)=\dfrac{\cos^3 \phi-\cos^5 \phi-\sin^4 \phi \cos\phi}{3}=\dfrac{\cos(3\phi)-\cos(5\phi)}{24}.
 \end{equation}
\end{lemma}
\proof
The assertion follows by expanding the $\sin$ and $\cos$ functions in the Taylor series about the origin and using~\eqref{eqE}.
\endproof

Now, for each fixed $\rho>0$, we represent $E(\rho,\phi)$ by its Fourier series
\begin{equation}\label{eqEFourier}
E(\rho,\phi)=a_0(\rho)+\sum\limits_{n=1}^\infty (a_n(\rho)\cos(n\phi)+b_n(\rho)\sin(n\phi))
\end{equation}
This series  converges to $E(\rho,\phi)$ for every $\rho$ and $\phi$ due to Lemma~\ref{lESmoothness}. Moreover, this lemma immediately implies the following result.

\begin{lemma}\label{lFourierE}
\begin{enumerate}
\item\label{lFourierE1} The Fourier coefficients $a_0(\rho)$, $a_n(\rho)$, and $b_n(\rho)$ and all their derivatives are continuous for $\rho>0$ and extend as continuous function to $\rho\ge 0$. Moreover, for any integer $M\ge 0$, there exists a constant $c_M>0$ such that
$$
 \left|a_n(\rho)\right|, \left|b_n(\rho)\right|,  \left|a_n'(\rho)\right|, \left|b_n'(\rho)\right|\le \dfrac{c_M}{n^M}\quad\text{for all } n\in\bbN,\ \rho\in(0,\pi].
$$
\item\label{lFourierE2} We have $a_3(0)=1/24$, $a_5(0)=-1/24$, and all the other Fourier coefficients of $E(0,\phi)$ vanish.
\end{enumerate}
\end{lemma}
In particular, Lemma~\ref{lFourierE} guarantees that the series in~\eqref{eqEFourier} converges absolutely and uniformly for $\rho\in[0,\pi]$ and $\phi\in[-\pi,\pi)$. In what follows, this will allow us to interchange the operations of summation and integration.

Now we are in a position to provide an asymptotics of $I_1(r,\psi)$.

\begin{lemma}\label{lI1}
   $I_1(r,\psi)= \dfrac{\pi\cos(4\psi)}{6 r^2} -\hat I_3(r,\psi)+ O(r^{-5/2})$ as $r\to\infty$, where $\hat I_3(r,\psi)$ is given by~\eqref{eqI3''}  and $O(\cdot)$ is uniform with respect to $\psi\in[-\pi,\pi)$.
\end{lemma}
\proof
Assume that $\psi\in\Psi_1$. Then, using integration by parts, we represent $I_1(r,\psi)$ as follows:
\begin{equation}\label{eqI1:1}
I_1(r,\psi)=\dfrac{1}{r\cos\psi}\int_{R_\pi}(\sin(x\theta))_{\theta_1}D_0(\theta)\,d\theta=
-\dfrac{1}{r\cos\psi}(\Sigma_1(r,\psi)+\Sigma_2(r,\psi)) - \hat I_3(r,\psi),
\end{equation}
where
$$
\Sigma_1(r,\psi):=\int_{B_\pi}\sin(x\theta)D(\theta)\,d\theta,\quad \Sigma_2(r,\psi):=\int_{R_\pi\setminus B_\pi}\sin(x\theta)\left(\dfrac{1}{A(\theta)}\right)_{\theta_1}\,d\theta,
$$
$D_0(\theta)$ and $D(\theta)$ are given by~\eqref{eqD0} and~\eqref{eqD}, respectively, and $\hat I_3(r,\psi)$ is given by~\eqref{eqI3''}.

Using~\eqref{eqEFourier}, we have
$$
\begin{aligned}
\Sigma_1(r,\psi)&=\int_0^\pi d\rho\int_{-\pi}^\pi \sin(r\rho\cos(\phi-\psi))E(\rho,\phi)\,d\phi\\
&=  \int_0^\pi \left(\sum\limits_{n=1}^\infty \int_{-\pi}^\pi \sin(r\rho\cos\phi)[a_n(\rho)\cos(n(\phi+\psi))+b_n(\rho)\sin(n(\phi+\psi))]\,d\phi\right)d\rho.
\end{aligned}
$$
Using the trigonometric addition formulas and the integral representation of the Bessel function $J_n$ (see~\eqref{eqBesselIntegralRepresentation2}), we obtain
$$
\begin{aligned}
\Sigma_1(r,\psi)& =  \int_0^\pi \left(\sum\limits_{\text{odd}\, n\in\bbN} \int_{-\pi}^\pi \sin(r\rho\cos\phi)\cos(n\phi)\,d\phi\right)[a_n(\rho)\cos(n\psi)+b_n(\rho)\sin(n\psi)]d\rho\\
&=2\pi\sum\limits_{\text{odd}\, n\in\bbN} (-1)^{(n-1)/2}   \left(\cos(n\psi)\int_0^\pi J_n(r\rho)a_n(\rho)\,d\rho+\sin(n\psi)\int_0^\pi J_n(r\rho) b_n(\rho)\,d\rho\right).
\end{aligned}
$$
Therefore, making the change of variables $z=r\rho$ in the integrals, taking into account Lemma~\ref{lFourierE} and applying Lemma~\ref{lSummationBesselLimitA}, we conclude that
\begin{equation}\label{eqI1:1'}
\begin{aligned}
\Sigma_1(r,\psi)& =\dfrac{2\pi}{r}\big( -\cos(3\psi)a_3(0)+\cos(5\psi)a_5(0) \big) + O(r^{-3/2})\\
&=-\dfrac{\pi}{12\, r}\big( \cos(3\psi)+\cos(5\psi)\big) + O(r^{-3/2}) \quad\text{as } r\to\infty,
\end{aligned}
\end{equation}
where $O(\cdot)$ is uniform with respect to $\psi$.

Now we estimate $\Sigma_2(r,\psi)$. Using the $2\pi$-periodicity of $\cos(x\theta)$ and $\frac{1}{A(\theta)}$ with respect to $\theta_1$, we obtain
\begin{equation}\label{eqI1:2}
\Sigma_2(r,\psi)=-\dfrac{1}{r\cos\psi}\int\limits_{R_\pi\setminus B_\pi}(\cos(x\theta))_{\theta_1}\left(\dfrac{1}{A(\theta)}\right)_{\theta_1}\,d\theta=
\tilde\Sigma_2(r,\psi)+\hat\Sigma_2(r,\psi),
\end{equation}
where
$$
\begin{aligned}
\tilde\Sigma_2(r,\psi)&:=\dfrac{\pi}{r \cos\psi} \int_{-\pi}^\pi  \cos(\pi r\cos(\phi-\psi))f(\phi)\cos\phi\,d\phi, \\
\hat\Sigma_2(r,\psi)&:=\dfrac{1}{r \cos\psi} \int\limits_{R_\pi\setminus B_\pi}\cos(x\theta)\left(\dfrac{1}{A(\theta)}\right)_{\theta_1\theta_1}\,d\theta,
\end{aligned}
$$
and $f(\phi)$ is obtained from $\left(\dfrac{1}{A(\theta)}\right)_{\theta_1}$ by substituting $\theta_1=\pi\cos\phi$, $\theta_2=\pi\sin\phi$.
Expanding $f(\phi)\cos\phi$ into the Fourier series
$$
f(\phi)\cos\phi=\alpha_0+\sum\limits_{n=1}^\infty \left(\alpha_n\cos(n\phi)+\beta_n\sin(n\phi)\right)
$$
and using the integral representation~\eqref{eqBesselIntegralRepresentation1},  we have
$$
\begin{aligned}
\tilde\Sigma_2(r,\psi)&=\dfrac{\pi}{r \cos\psi} \int_{-\pi}^\pi  \cos(\pi r\cos\phi)\alpha_0 d\phi\\
 & + \dfrac{\pi}{r \cos\psi} \sum\limits_{\text{even}\, n=2}^\infty  \int_{-\pi}^\pi  \cos(\pi r\cos\phi)\left(\alpha_n\cos(n\phi)\cos(n\psi) + \beta_n \cos(n\phi)\sin(n\psi)\right)  d\phi\\
  &= \dfrac{\pi^2}{r \cos\psi} \left(\alpha_0 J_0(\pi r) +  (-1)^{n/2}\sum\limits_{\text{even}\, n=2}^\infty \left(\alpha_n \cos(n\psi)J_n(\pi r) + \beta_n \sin(n\psi)J_n(\pi r)\right)\right)
\end{aligned}
$$
Thus, Lemma~\ref{lSummationBesselLimitAlpha} implies
\begin{equation}\label{eqI1:3}
\tilde\Sigma_2(r,\psi)=O(r^{-3/2}),
\end{equation}
where   $O(\cdot)$ is uniform with respect to $\psi\in[-\pi,\pi)$.

Finally, using integration by parts, we immediately obtain $\hat\Sigma_2(r,\psi)=O(r^{-2}).$ Combining this with~\eqref{eqI1:1}--\eqref{eqI1:3} yields
$$
\begin{aligned}
I_1(r,\psi)&=\dfrac{\pi}{12\, r^2\cos\psi}\big( \cos(3\psi)+\cos(5\psi)\big) -\hat I_3(r,\psi)+ O(r^{-5/2})\\
&=\dfrac{\pi\cos 4\psi}{6 r^2} -\hat I_3(r,\psi)+ O(r^{-5/2}).
\end{aligned}
$$
\endproof

Now Theorem~\ref{thAsympfMultiDim} follows from~\eqref{eqOmegaI1I3I4} and Lemmas~\ref{lI3}, \ref{lI4}, and~\ref{lI1}.

\appendix

\section{Auxiliary results}\label{appendixAuxiliary}

In this appendix, we collect some known facts about  the Bessel functions of the first kind $J_n(z)$ as well as several corollaries that we need in the main part of the paper. We will consider the functions $J_n(z)$  only for $z>0$ and integer $n\ge 0$.

\subsection{Known facts}

\begin{lemma}[see~{\cite[Sec.~9.1.60]{Abramowitz}} and~{\cite[Sec.~10.22.41]{Olver}}]
For all integer $n\ge 0$, we have
\begin{align}
&|J_n(z)|\le 1\quad\text{for all } z\ge 0, \label{eqBesselProperty1}\\
& \int_0^\infty J_n(z)\,dz = 1. \label{eqBesselProperty2}
\end{align}
\end{lemma}

Let $\gamma$ denote the Euler--Mascheroni constant.

\begin{lemma}[see {\cite[Sec. 12.2, Example 4]{Whittaker} and~\cite[Section 11.1.20]{Abramowitz}}]\label{lExpGamma} We have
  $$ \int_0^1\dfrac{1-e^{-z}}{z}\,dz -
 \int_1^\infty\dfrac{e^{-z}}{z}\,dz =   \int_0^1
\dfrac{1-J_0(z)}{z}\,dz -  \int_1^{\infty}
\dfrac{J_0(z)}{z}\,dz + \ln 2= \gamma.
 $$
\end{lemma}

The following lemma provides two integral representations of the Bessel functions.

\begin{lemma}[see {\cite[Section 9.1.21]{Abramowitz}}]\label{lBesselIntegralRepresentation}
 We have
    \begin{align}
  & J_n(z)=\dfrac{(-1)^{n/2}}{\pi}\int_0^\pi \cos(z\cos\phi)\cos(n\phi)\,d\phi & &\text{for even $n\ge 0$},\label{eqBesselIntegralRepresentation1}\\
  & J_n(z)=\dfrac{(-1)^{(n-1)/2}}{\pi}\int_0^\pi \sin (z\cos\phi)\cos(n\phi)\,d\phi & &\text{for odd $n\ge 1$}.\label{eqBesselIntegralRepresentation2}
  \end{align}

\end{lemma}

We will use the following asymptotics for the function $J_0(z)$.

\begin{lemma}[see~{\cite[Section~7.21, p. 199]{Watson}}]\label{lAsympBessel0} We have
  $$
J_0(z)=\dfrac{\sqrt{2}}{\sqrt{\pi} z^{1/2}} \cos\left(z-\dfrac{\pi}{4}\right)+
\dfrac{\sqrt{2}}{8\sqrt{\pi} z^{3/2}}\sin\left(z-\dfrac{\pi}{4}\right)+ O\left(\dfrac{1}{z^{5/2}}\right)
$$
as $z\to\infty$.
\end{lemma}

A similar asymptotics holds for $J_n(z)$. We will need it in the following form
\begin{lemma}[see~{\cite[Section~7.21, p. 199]{Watson}}]\label{lAsympBesseln} For each $n\ge 0$, there exist a constant $C_n>0$ such that, for all $z\ge 1$,
\begin{align}
 & J_n(z)=\left(\dfrac{2}{\pi z}\right)^{1/2}\cos\left(z-\frac{\pi n}{2}-\frac{\pi}{4}\right) + R_n(z), \label{eqBesselAsympn'}\\
 & |R_n(z)|\le \dfrac{C_n}{z^{3/2}}. \label{eqBesselAsympn''}
\end{align}
\end{lemma}

\subsection{Corollaries from the known facts}
The following result provides an explicit dependence of the remainder $R_n(z)$ in~\eqref{eqBesselAsympn'} for large $z$, namely $z>n^2$.

\begin{lemma}\label{lBesselAsympGeneral} There is a constant $C>0$ such that  for all integer $n\ge 0$ and for all $z>n^2$, we have
\begin{align}
 & J_n(z)=\left(\dfrac{2}{\pi z}\right)^{1/2}\cos\left(z-\frac{\pi n}{2}-\frac{\pi}{4}\right) + R_{n}(z), \label{eqBesselAsympGeneral'}\\
 & |R_{n}(z)|\le C n^3\dfrac{1}{z^{3/2}}. \label{eqBesselAsympGeneral''}
\end{align}
\end{lemma}
\proof
Consider the Hankel functions $H_n^{1}(z)$ and $H_n^{2}(z)$ that are connected with the Bessel functions $J_n(z)$ via
\begin{equation}\label{eqHankelBessel}
  J_n(z)=\dfrac{H_n^{1}(z)+H_n^{2}(z)}{2}
\end{equation}
and have the following asymptotics (see~\cite[Section~7.2]{Watson}):
\begin{equation}\label{eqBesselAsympGeneral1}
\begin{aligned}
H_n^{(1)}(z)&=\left(\dfrac{2}{\pi z}\right)^{1/2}e^{i\left(z-\frac{\pi n}{2}-\frac{\pi}{4}\right)}\left[\sum\limits_{m=0}^{p-1}\dfrac{(\frac{1}{2}-n)_m \, \Gamma\left(n+m+\frac{1}{2}\right)}{m!\,\Gamma\left(n+\frac{1}{2}\right)(2iz)^m} + R_{n,p}^{(1)}(z)\right],\\
H_n^{(2)}(z)&=\left(\dfrac{2}{\pi z}\right)^{1/2}e^{-i\left(z-\frac{\pi n}{2}-\frac{\pi}{4}\right)} \left[\sum\limits_{m=0}^{p-1}\dfrac{(\frac{1}{2}-n)_m \, \Gamma\left(n+m+\frac{1}{2}\right)}{m!\,\Gamma\left(n+\frac{1}{2}\right)(-2iz)^m} + R_{n,p}^{(2)}(z)\right].
\end{aligned}
\end{equation}
For an arbitrary fixed constant $\sigma\in(0,1)$ and $p\ge n$, the remainder $R_{n,p}^{(1)}(z)$ is estimated for all $z>0$ as follows:
$$
|R_{n,p}^{(1)}(z)|\le \dfrac{\sigma^{n-p-1/2}\left|\left(\frac{1}{2}-n\right)_p\right|}{(p-1)!\,\Gamma\left(n+\frac{1}{2}\right)(2z)^p}\int_0^\infty e^{-u}u^{n+p-1/2}du.
$$
Therefore, taking $p=n$ and $z> n^2$, we obtain
\begin{equation}\label{eqBesselAsympGeneral2}
\begin{aligned}
|R_{n,n}^{(1)}(z)|&\le \dfrac{\sigma^{-1/2}\left|\left(\frac{1}{2}-n\right)\left(\frac{1}{2}-n+1\right)\cdot{\dots}\cdot\left(\frac{1}{2}-1\right)\right|\,\Gamma\left(2n+\frac{1}{2}\right)}
{(n-1)!\,\Gamma\left(n+\frac{1}{2}\right)2^n  n^{2n-2} z}\,\\
&\le \dfrac{\sigma^{-1/2}\,(2n)!}
{(n-1)!\,2^n  n^{2n-3} z}
\le \dfrac{\sigma^{-1/2}\,n(n+1)\cdot{\dots}\cdot(2n)}
{2^n  n^{2n-3} z}\le \dfrac{\sigma^{-1/2}\,(2n)^{n+1}}
{2^n  n^{2n-3} z}\le \dfrac{2\sigma^{-1/2}}
{n^{n-4} z}.
\end{aligned}
\end{equation}

Further, for $z>n^2$, we have
\begin{equation}\label{eqBesselAsympGeneral3}
\begin{aligned}
 & \left|e^{i\left(z-\frac{\pi n}{2}-\frac{\pi}{4}\right)}\sum\limits_{m=1}^{n-1}\dfrac{(\frac{1}{2}-n)_m \, \Gamma\left(n+m+\frac{1}{2}\right)}{m!\,\Gamma\left(n+\frac{1}{2}\right)(2iz)^m}\right| \le \sum\limits_{m=1}^{n-1}\dfrac{n\cdot{\dots}\cdot (n-m+1) \, (n+m)!}{m!\,(n-1)!\,2^m n^{2m-2}z}\\
&= \sum\limits_{m=1}^{n-1}\dfrac{(n+m)!}{m!\,(n-m)!\,2^m n^{2m-3}z}
= \sum\limits_{m=1}^{n-1}\dfrac{\binom{n}{m} (n+1)\cdot{\dots}\cdot(n+m)}{2^m n^{2m-3}z}\\
&\le  \sum\limits_{m=1}^{n-1}\dfrac{\binom{n}{m}}{n^{m-3}z}\le n^3 \left(1+\frac{1}{n}\right)^n \dfrac{1}{z}.
\end{aligned}
\end{equation}
Combining~\eqref{eqBesselAsympGeneral1}--\eqref{eqBesselAsympGeneral3}, we obtain for all integer $n\ge 0$ and $z\ge n^2$
\begin{equation}\label{eqBesselAsympGeneral4}
\begin{aligned}
 & H_n^{(1)}(z)=\left(\dfrac{2}{\pi z}\right)^{1/2}e^{i\left(z-\frac{\pi n}{2}-\frac{\pi}{4}\right)}\left[1+S_{n}^{(1)}(z)\right],\\
 & |S_{n}^{(1)}(z)|\le C n^3\dfrac{1}{z},
\end{aligned}
\end{equation}
where $C>0$ does not depend on  $z\ge n^2$.

Combining~\eqref{eqBesselAsympGeneral4} with the analogous representation for $H_n^{(2)}$ and taking into account~\eqref{eqHankelBessel}, we complete the proof.
\endproof

In the following lemmas, we deal with series involving the Bessel functions $J_n(z)$ and the sequences $\alpha_n$ and $A_n(z)$. In the main part of the text, these sequences arise as the Fourier coefficients of certain smooth  functions.

\begin{lemma}\label{lSummationBesselLimitAlpha}
    Let a sequence $\alpha_n\in\bbR$, $n\in\bbN$, satisfy
\begin{equation}\label{eqSummationBesselLimit'}
   |\alpha_n|\le \dfrac{c_0}{n^5},
\end{equation}
where $c_0>0$ does not depend on $n$. Then
\begin{equation}\label{eqSummationBesselLimit}
\sum\limits_{n=1}^\infty \alpha_n J_n(z)=\left(\dfrac{2}{\pi z}\right)^{1/2}\sum\limits_{n=1}^\infty \alpha_n \cos\left(z-\frac{\pi n}{2}-\frac{\pi}{4}\right) + O\left(\frac{1}{z^{3/2}}\right)\quad\text{as } z\to\infty.
\end{equation}
\end{lemma}
\proof Let $R_n(z)$ be the remainder in the asymptotics for $J_n(z)$ defined in Lemmas~$\ref{lAsympBesseln}$ and~$\ref{lBesselAsympGeneral}$. Then we have to prove that
\begin{equation}\label{eqSummationBesselLimitAlpha1}
\sum\limits_{n=1}^\infty \alpha_n z^{3/2} R_n(z) = O(1) \quad\text{as } z\to\infty.
\end{equation}

Let $1<z<n^2$. Then, using~\eqref{eqBesselProperty1} and~\eqref{eqSummationBesselLimit'}, we have
\begin{equation}\label{eqSummationBesselLimitAlpha2}
\begin{aligned}
\left| z^{3/2} R_n(z)\alpha_n\right| & \le n^3 \left|J_n(z)-\left(\dfrac{2}{\pi z}\right)^{1/2}\cos\left(z-\frac{\pi n}{2}-\frac{\pi}{4}\right)\right| |\alpha_n|\\
& \le n^3\left(1+\left(\dfrac{2}{\pi}\right)^{1/2}\right) \dfrac{c_0}{n^5} \le \left(1+\left(\dfrac{2}{\pi}\right)^{1/2}\right) \dfrac{c_0}{n^2}.
\end{aligned}
\end{equation}

Let $z\ge n^2$. Then, using Lemma~\ref{lBesselAsympGeneral} and inequalities~\eqref{eqSummationBesselLimit'}, we obtain
\begin{equation}\label{eqSummationBesselLimitAlpha3}
\left| z^{3/2} R_n(z)\alpha_n\right|   \le z^{3/2} Cn^3\dfrac{1}{z^{3/2}} \dfrac{c_0}{n^5}= \dfrac{C c_0}{n^2}.
\end{equation}

Estimates~\eqref{eqSummationBesselLimitAlpha2} and~\eqref{eqSummationBesselLimitAlpha3} imply~\eqref{eqSummationBesselLimitAlpha1}.
\endproof

\begin{lemma}\label{lSummationBesselLimitA}   Assume a sequence $A_n\in C^1[0,1]$, $n\in\bbN$, satisfies
\begin{equation}\label{eqSummationBessel1}
  \sup\limits_{\rho\in[0,1]} |A_n(\rho)|\le \dfrac{c_0}{n^2},\quad \sup\limits_{\rho\in[0,1]}|A_n'(\rho)|\le \dfrac{c_1}{n^6} \quad\text{for all } n\in\bbN,
\end{equation}
where $c_0,c_1>0$ do not depend on $n$. 
Furthermore, assume that there exists a finite set $\cN\subset\bbN$ such that
\begin{equation}\label{eqSummationBessel1'}
A_n(0)=0\quad\text{for all } n\notin\cN.
\end{equation}
Then there is a constant $c>0$ depending only on $\cN$ and on the constants $c_0$ and $c_1$ such that
$$
\left|\int_0^{\sigma} \left(\sum\limits_{n=1}^\infty  J_n(z)A_n\left(\frac{z}{\sigma}\right)\right)dz - \sum\limits_{n\in\cN} A_n(0)\right| \le \dfrac{c}{\sigma^{1/2}}\quad \text{for all } \sigma>1.
$$
\end{lemma}
\proof
For any fixed $\sigma>1$, inequalities~\eqref{eqBesselProperty1} and~\eqref{eqSummationBessel1} imply that the series $\sum\limits_{n=1}^\infty  J_n(z)A_n\left(\frac{z}{\sigma}\right)$ converges uniformly with respect to $z\in(0,\sigma)$. Hence,
\begin{equation}\label{eqSummationBesselLimit1'}
\begin{aligned}
& \int_0^{\sigma} \left(\sum\limits_{n=1}^\infty  J_n(z)A_n\left(\frac{z}{\sigma}\right)\right)dz  =\sum\limits_{n=1}^\infty  \int_0^{\sigma} J_n(z)A_n\left(\frac{z}{\sigma}\right) dz\\
&\qquad = \sum\limits_{n\in\cN}  \int_0^{\sigma} J_n(z)A_n(0)\,dz + \sum\limits_{n=1}^\infty  \int_0^{\sigma} J_n(z)\left(A_n\left(\frac{z}{\sigma}\right)-A_n(0)\right) dz.
\end{aligned}
\end{equation}
Consider the first sum in the right-hand side in~\eqref{eqSummationBesselLimit1'}.
Due to~\eqref{eqBesselProperty2},
$$
  \int_0^{\sigma} J_n(z) dz=1-\int_\sigma^\infty J_n(z) dz.
$$
Using Lemma~\ref{eqBesselAsympn'} and integration by parts, we obtain
$$
\left|\int_\sigma^\infty J_n(z) dz\right|\le \left|\int_\sigma^\infty \left(\dfrac{2}{\pi z}\right)^{1/2}\cos\left(z-\frac{\pi n}{2}-\frac{\pi}{4}\right)dz\right| + \int_\sigma^\infty \dfrac{C_n}{z^{3/2}}\,dz \le \dfrac{\tilde C_n}{\sigma^{1/2}},
$$
where $\tilde C_n>0$ does not depend on $\sigma>0$. Therefore,
\begin{equation}\label{eqSummationBesselLimit1}
  \left|\sum\limits_{n\in\cN}\int_0^{\sigma} J_n(z)A_n(0) dz-\sum\limits_{n\in\cN}A_n(0)\right|\le   \dfrac{c_2c_0}{\sigma^{1/2}},
\end{equation}
where $c_2>0$ depends on $\tilde C_n$, $n\in\cN$, but does not depend on $\sigma$.

Consider the second sum the right-hand side in~\eqref{eqSummationBesselLimit1'}. Using~\eqref{eqBesselAsympn'}, consider for $\sigma>1$
\begin{equation}\label{eqSummationBesselLimit2}
  \int_0^{\sigma} J_n(z)\left(A_n\left(\frac{z}{\sigma}\right)-A_n(0)\right) dz=G_1(\sigma) + \left(\dfrac{1}{\pi}\right)^{1/2} (G_2(\sigma)+G_3(\sigma)),
\end{equation}
where
$$
\begin{aligned}
G_1(\sigma)&:=\int_0^1 J_n(z)\left(A_n\left(\frac{z}{\sigma}\right)-A_n(0)\right) dz,\\
G_2(\sigma)&:=\int_1^{\sigma} \dfrac{1}{z^{1/2}}\cos\left(z-\frac{\pi n}{2}-\frac{\pi}{4}\right) \left(A_n\left(\frac{z}{\sigma}\right)-A_n(0)\right) dz,\\
G_3(z)&:=\int_1^{\sigma}  R_n(z) \left(A_n\left(\frac{z}{\sigma}\right)-A_n(0)\right)dz.
\end{aligned}
$$
Due to~\eqref{eqBesselProperty1}, we have
\begin{equation}\label{eqSummationBesselLimit3}
  |G_1(\sigma)|\le \sup\limits_{\rho\in[0,1]} |A_n'(\rho)|\dfrac{1}{\sigma}\le\dfrac{c_1}{n^2}\dfrac{1}{\sigma}.
\end{equation}
To estimate $G_2(\sigma)$, we integrate by parts and use~\eqref{eqSummationBessel1}:
\begin{equation}\label{eqSummationBesselLimit4}
\begin{aligned}
  |G_2(\sigma)| &  \le \left|\dfrac{1}{z^{1/2}}\sin\left(z-\frac{\pi n}{2}-\frac{\pi}{4}\right)(A_n\left(\frac{z}{\sigma}\right)-A_n(0))\Big|_{z=1}^{z=\sigma}\right| \\
& +
\int_{1}^{\sigma}\dfrac{1}{2z^{3/2}}\left|(A_n\left(\frac{z}{\sigma}\right)-A_n(0)\right|\,dz + \dfrac{1}{\sigma}\int_{1}^{\sigma}\dfrac{1}{z^{1/2}}\left|A_n'\left(\frac{z}{\sigma}\right)\right|\,dz\\
&\le  \dfrac{|A_n(1)-A_n(0)|}{\sigma^{1/2}}  + \left|A_n\left(\frac{1}{\sigma}\right)-A_n(0)\right| \\
 & + \sup\limits_{\rho\in[0,1]} |A_n'(\rho)| \int_{1}^{\sigma}\dfrac{1}{2z^{3/2}}\dfrac{z}{\sigma}\,dz +  \dfrac{c_3}{\sigma^{1/2}} \sup\limits_{\rho\in[0,1]} |A_n'(\rho)|\le \dfrac{c_4}{n^2\sigma^{1/2}},
\end{aligned}
\end{equation}
where $c_3,c_4>0$ do not depend on $n$ or $\sigma$.

Now let us estimate $G_3(\sigma)$. If $\sigma<n^2$, then, using the inequalities (see~\eqref{eqBesselProperty1})
$$
|R_n(z)|\le |J_n(z)|+\left(\frac{2}{\pi}\right)^{1/2}\le 1+\left(\frac{2}{\pi}\right)^{1/2}=:c_5\quad\text{for all } z>1,
$$
and~\eqref{eqSummationBessel1}, we have
\begin{equation}\label{eqSummationBesselLimit5}
|G_3(\sigma)| \le c_5\sup\limits_{\rho\in[0,1]}|A_n'(\rho)|\int_1^{n^2} \dfrac{z}{\sigma}\,dz \le c_5 \dfrac{c_1}{n^6}\dfrac{n^4}{2\sigma}\le
 \dfrac{c_5 c_1}{2n^2}\dfrac{1}{\sigma}.
\end{equation}
If $\sigma>n^2$, then we write
\begin{equation}\label{eqSummationBesselLimit6}
G_3(\sigma) = G_{3,1}(\sigma)+G_{3,2}(\sigma),
\end{equation}
where
$$
G_{3,1}(\sigma):=\int_1^{n^2}  R_n(z) \left(A_n\left(\frac{z}{\sigma}\right)-A_n(0)\right)dz,\quad  G_{3,2}(\sigma):=\int_{n^2}^{\sigma}  R_n(z) \left(A_n\left(\frac{z}{\sigma}\right)-A_n(0)\right)dz.
$$
Similarly to~\eqref{eqSummationBesselLimit5}, we have
\begin{equation}\label{eqSummationBesselLimit7}
|G_{3,1}(\sigma)|\le  \dfrac{c_5 c_1}{2n^2}\dfrac{1}{\sigma}.
\end{equation}
Finally, using~\eqref{eqSummationBessel1} and Lemma~\ref{lBesselAsympGeneral}, we obtain
\begin{equation}\label{eqSummationBesselLimit8}
|G_{3,2}(\sigma)|\le \sup\limits_{\rho\in[0,1]}|A_n'(\rho)|\int_{n^2}^\sigma \dfrac{z}{\sigma} Cn^3\dfrac{1}{z^{3/2}}\,dz\le
\dfrac{c_1}{n^6} \dfrac{Cn^3}{\sigma}2\sigma^{1/2}\le \dfrac{2c_1C}{n^3}\dfrac{1}{\sigma^{1/2}}.
\end{equation}

Combining~\eqref{eqSummationBesselLimit2}--\eqref{eqSummationBesselLimit8}, we see that
$$
  \sum\limits_{n=1}^\infty\left|\int_0^{\sigma} J_n(z)\left(A_n\left(\frac{z}{\sigma}\right)-A_n(0)\right) dz\right|\le \dfrac{c_6(c_0+c_1)}{\sigma^{1/2}}.
$$
where $c_6>0$ does not depend on $\sigma$.
This relation, together with~\eqref{eqSummationBesselLimit1'} and~\eqref{eqSummationBesselLimit1} completes the proof.
\endproof

{\bf Acknowledgement.} The author expresses his gratitude to
Sergey Tikhomirov for numerous discussions and to an anonymous referee for meticulously reading the manuscript and providing useful suggestions. The research was
supported by the DFG project SFB 910, the DFG Heisenberg Programme and by the ``RUDN University Program 5-100''.

\end{document}